\documentclass[12pt,reqno]{amsart}

\usepackage{ytableau}
\usepackage{slashbox}
\usepackage{fullpage}
\usepackage{epsfig}
\usepackage{amsfonts,amssymb,amscd,amsmath,euscript,enumerate,verbatim,calc}
\usepackage{chngcntr}
\theoremstyle{plain}% default
\newtheorem{theorem}{Theorem}

\newtheorem{lemma}{Lemma}

\theoremstyle{definition}
\newtheorem{definition}{Definition}
\newtheorem{example}{Example}

\newtheorem{conjecture}{Conjecture}

\newenvironment{prf}[1][Proof]{\noindent \textbf{#1.} }{\  \hfill \rule{0.5em}{0.5em}}
\newdimen \dummy
\dummy=\oddsidemargin
\newcommand{\rmv}[1]{}
\counterwithin{theorem}{section}
\counterwithin{lemma}{section}
\counterwithin{definition}{section}
\counterwithin{example}{section}

\begin{document}

\title{Core Partitions With d-Distinct Parts}

\author{Murat Sahin}

\address{Department of Mathematics,  Faculty of Sciences,
          Ankara University, Tandogan, Ankara, 06100, Turkey}
\email{msahin@ankara.edu.tr}

\date{}



\begin{abstract}
In this paper, we study $(s,s+1)$-core partitions with $d$-distinct parts.
We obtain results on the number and the largest size of such partitions,
so we extend Xiong's paper in which the results are obtained about
$(s,s+1)$-core partitions with distinct parts. Also, we propose a conjecture about
$(s,s+r)$-core partitions with $d$-distinct parts for $1 \le r \le d$.
\rmv{we propose a conjecture about $(s,ts+r)$-core partitions with $d$-distinct
parts for $1 \le r \le d$.}
\end{abstract}

\maketitle

\section{Introduction}

A \emph{partition} of $n$ is a finite nonincreasing sequence
$\lambda=(\lambda_1,\lambda_2,\ldots, \lambda_l)$ such that
$n=\lambda_1+\lambda_2+ \ldots +\lambda_l$. A summand in a
partition is called \emph{part}. We say that $n$ is the size of
$\lambda$ and $l$ is the length of $\lambda$. For example,
$\lambda=(6, 3, 3, 2, 1)$ is a partition of $n=15$. The parts of
the partition $\lambda$ are $6, 3, 3, 2$ and $1$. The size of $\lambda$
is $15$ and the length of $\lambda$ is $5$. A partition $\lambda$ is
called \emph{a partition with $d$-distinct parts} if and only if
$\lambda_i-\lambda_{i+1} \ge d$.

Partitions can be visualized with \emph{Young diagram}, which is a
finite collection of boxes arranged in left-justified rows, with $\lambda_i$ boxes
in the $i$-th row. The pair $(i,j)$ shows that the coordinates of the boxes in the Young diagram.
The Young diagram of  $\lambda=(6, 3, 3, 2, 1)$
is as follows:

\begin{equation*}\label{diagram}
\ytableausetup{boxsize=2.3em}
\begin{ytableau}
\tiny{(1,1)} & (1,2)& (1,3)& (1,4) &(1,5) & (1,6) \cr (2,1) & (2,2) & (2,3) \cr (3,1) & (3,2)& (3,3) \cr (4,1) & (4,2) \cr (5,1) \cr
\end{ytableau}
\end{equation*}

For each boxes in the Young diagram in coordinates $(i,j)$, the \emph{hook length} is defined as
the sum of the number of boxes exactly to the right, exactly
below and the box itself. So, the hook lengths of  the partition $\lambda=(6, 3, 3, 2, 1)$
can be given as follows:

\begin{equation} \label{hookdiagram}
\begin{ytableau}
10 & 8 & 6 & 3 & 2 & 1 \\
6 & 4 & 2  \\
5 & 2 & 1  \\
2 & 1      \\
1     \\
\end{ytableau}
\end{equation}

$h(i,j)$ will show that the entry in the coordinate $(i,j)$ of the box, that is
the hook length of the box. If $\lambda=(6, 3, 3, 2, 1)$, then $h(1,1)=10, h(2,3)=2$ and
$h(6,1)=1$ as you can see from~(\ref{hookdiagram}).

A partition $\lambda$ is called a
\emph{$s$-core} partition if $\lambda$ has no boxes of hook
length $s$. For example, the partition $\lambda=(6, 3, 3, 2, 1)$
is $7$-core but it is not $5$-core, since  $\lambda$ has no boxes
of hook length $7$, but it has a box of hook length $5$ (See the diagram~(\ref{hookdiagram})).

A more general definition, a partition
$\lambda$ is called a \emph{($s_1,s_2,\ldots,s_t$)-core}
partition if $\lambda$ has no boxes of hook length $s_1,s_2,\ldots,s_t$.
So, the partition $\lambda=(6, 3, 3, 2, 1)$ is $(7,9)$-core.

There are many studies about core partitions and such partitions
are closely related to posets, cranks, Raney numbers,
Catalan numbers, Fibonacci numbers etc. \cite{amdeberhan, cranks, posets, raney}.

Anderson \cite{uc} gives the result that the number of $(s,t)$-core partitions is finite
if and only if $s$ and $t$ are coprime. In that case, this number is $$\frac{1}{s+t}\binom{s+t}{s}.$$
Olsson and Stanton \cite{OS} give the largest size of such partitions. Some results on the number,
the largest size and the average size of such partitions are provided by \cite{AmdeLeven,Armstrong,Johnson,NathSellers,StanleyZanello,Xiong1,Xiong2,posets}.
In particular, the number of $(s,s+1)$-core partitions is the Catalan number $$C_s=\frac{1}{s+1}\binom{2s}{s}.$$

Amdeberhan \cite{amdeberhan} conjectures that the number of $(s,s+1)$-core partitions with
distinct parts equals the Fibonacci numbers. This conjecture is proved
independently by Xiong \cite{xiong} and Straub \cite{straub}. More generally, Straub \cite{straub} characterize 
the number $N_d(s)$ of $(s,ds-1)$-core partitions with distinct parts by $N_d(1)=1$,$N_d(2)=d$ and, for $s \geq 3$,
$$N_d(s)=N_d(s-1)+dN_d(s-2).$$
Xiong \cite{xiong} also obtain results
on the number, the largest size and the avarage size of $(s,s+1)$-core partitions with
distinct parts. 

In this paper, we ask the problem of counting the number of special partitions which are $s$-core for certain values of $s$.
More precisely, we focus on $(s,s+1)$-core partitions with $d$-distinct parts.
We obtain results on the number, the largest size of such partitions,
so we extend Xiong's paper in which the results are obtained about $(s,s+1)$-core partitions
with distinct parts. Also, we propose the following conjecture about
$(s,s+r)$-core partitions with $d$-distinct parts for $1 \le r \le d$. That is, the number $N_{d,r}(s)$ of $(s,s+r)$-core partitions with
$d$-distinct parts is characterized by
$N_{d,r}(s)=s$ for $1 \le s \le d$, $N_{d,r}(d+1)=d+r$, and for $s \ge d+2$,
$$ N_{d,r}(s)=N_{d,r}(s-1)+N_{d,r}(s-(d+1))$$
for $1 \le r \le d$.

\section{$(s,s+1)$-core partitions with $d$-distinct parts}
Suppose that $\lambda=(\lambda_1,\lambda_2,\ldots, \lambda_l)$ is a
partition whose corresponding Young diagram has $l$ rows. The set
$\beta(\lambda)$ of $\lambda$ is defined to be the set of first
column hook length in the Young diagram of $\lambda$, i.e.,
$$ \beta(\lambda)= \left\lbrace h(i,1):1 \le i \le l \right\rbrace. $$
For example, if $\lambda=(6, 3, 3, 2, 1)$, then we get
\begin{eqnarray*}
\beta(\lambda) & = & \left\lbrace h(1,1), h(2,1), h(3,1), h(4,1), h(5,1) \right\rbrace \\
               & = & \left\lbrace 10, 6, 5, 2, 1 \right\rbrace
\end{eqnarray*}
by using the diagram~(\ref{hookdiagram}).

Now, we generalize the definition of the twin-free set in \cite{straub}.
\begin{definition} \label{freeset}
Suppose that $d$ is a positive integer such that $d \geq 2$. A set
$X \subseteq \mathbb{N}$ is called $d$-th order twin-free set if there is no $x \in X$ such that
$$\left\lbrace x,x+k\right\rbrace \subseteq X, \quad \text{for} \:\: 1 \leq k \leq d. $$
\end{definition}

If we take $d=1$ in Definition~\ref{freeset}, then we get the twin-free set in \cite{straub}, that is,
we called it first order twin-free set.

\begin{example}
Let us take	$X=\left\lbrace 10, 5,2 \right\rbrace$. For $d=2$, $X$ is a second order twin-free set, since the sets
$$
\left\lbrace 2, 3\right\rbrace, \left\lbrace 5, 6\right\rbrace,
\left\lbrace 10, 11\right\rbrace,\left\lbrace 2, 4\right\rbrace,
\left\lbrace 5, 7\right\rbrace, \left\lbrace 10, 12\right\rbrace
$$
are not a subset of $X$. But, the set $\left\lbrace 2, 5\right\rbrace$ is a subset of $X$,
so $X$ is not a third order twin-free set.
 
\end{example}

\begin{theorem}\label{Abacus}
	\begin{enumerate}[(i)]
		\item Suppose $\lambda=(\lambda_1,\lambda_2,\ldots, \lambda_l)$ is a partition. Then
		$$\lambda_i=h(i,1)-l+1, \quad \text{for} \:\: 1 \le i \le l.$$
		Thus, $$|\lambda|=\sum_{x \in \beta(\lambda)}x-\binom{l}{2}.$$
		\item A partition $\lambda$ is a $s$-core partition if and only if for any $x \in \beta(\lambda)$
		with $x > s$, we always have $x-s \in \beta(\lambda)$.    
	\end{enumerate}
\end{theorem}

\begin{prf}
	See \cite{uc,bes}.
\end{prf}

\begin{lemma}\label{free}
The partition $\lambda$ is a partition with $d$-distinct parts if and only if $\beta(\lambda)$ is a $d$-th order twin free set. 	
\end{lemma}

\begin{prf}
Suppose that $\lambda=(\lambda_1,\lambda_2,\ldots, \lambda_l)$ is a partition.
If $\lambda$ is a partition with $d$-distinct parts if and only if
$\lambda_i - \lambda_{i+1} \ge d$. Then by Theorem~\ref{Abacus}(i),
\begin{eqnarray*}
h(i,1)-h(i+1,1)&=& (\lambda_i +l-1)-(\lambda_{i+1} +l-1) \\
           &=& \lambda_i-\lambda_{i+1} \\
           &\ge& d.
\end{eqnarray*}
So,  we get
\begin{eqnarray*}
	h(i,1)-h(i+1,1) \ge d & \iff & \beta(\lambda) \: \text{is a} \: d\text{-th order twin-free set}. \\
\end{eqnarray*}
\end{prf}

\begin{lemma} \label{bound}
	Suppose $\lambda$ is a $(s,s+1)$-core partition with $d$-distinct parts. Then, $$\beta(\lambda) \subset \{1,2,\ldots, s-1\}.$$
\end{lemma}

\begin{prf}
Suppose that $\lambda$ is a partition with $d$ distinct parts.
Then, $\beta(\lambda)$ is a $d$-th order twin-free set by Lemma~\ref{free}.
Since $\lambda$ is a $(s,s+1)$-core partition, we have
$ s,s+1\notin \beta(\lambda)$. If $x \ge s+2$ and $x \in \beta(\lambda)$,
then by Theorem~\ref{Abacus}(ii), we know that $x-s, x-(s+1) \in \beta(\lambda)$.
But, this is a contradiction since $\beta(\lambda)$ is $d$-th order twin-free set. That is,
$x \notin \beta(\lambda)$ and so we get the required result   
$$\beta(\lambda) \subset \{1,2,\ldots, s-1\}.$$
\end{prf}

\begin{lemma}\label{lemma}
	A partition $\lambda$ is a $(s,s+1)$-core partition with
	$d$-distinct parts if and only if $\beta(\lambda)$ is
	a $d$-th order twin-free subset of the set $\{1,2,\ldots, s-1\}$.	
\end{lemma}

\begin{prf}
If a partition $\lambda$ is a $(s,s+1)$-core partition with
$d$-distinct parts then by Lemma~\ref{bound}, $\beta(\lambda)$
must be subset of the set $\{1,2,\ldots, s-1\}$. Also, By
Lemma~\ref{free}, $\beta(\lambda)$ must be a $d$-th order
twin-free set.

Conversely, suppose that $\beta(\lambda)$ is a $d$-th order twin-free
subset of the set $\{1,2,\ldots, s-1\}$. By Lemma~\ref{free},
$\lambda$ is a partition with $d$-distinct parts. Also,
since $\beta(\lambda)$ is a subset of the set $\{1,2,\ldots, s-1\}$,
all the hook lengths of corresponding partition smaller than $s$ and $s+1$.
This means that $\lambda$ is a $(s,s+1)$-core partition. 
	
\end{prf}

\begin{theorem}\label{number}
The number $N_{d}(s)$ of $(s,s+1)$-core partitions with $d$-distinct parts is characterized by
	$N_{d}(s)=s$ for $1 \le s \le (d+1)$, and for $s \ge d+2$,
	$$ N_{d}(s)=N_{d}(s-1)+N_{d}(s-(d+1)).$$
\end{theorem}

\begin{prf}
	
Let $X_k$ denote the set of all $d$-th order twin-free
subsets of the set $\{1,2,\ldots, k-1\}$.
A partition $\lambda$ is a $(s,s+1)$-core partition with
$d$-distinct parts if and only if $\beta(\lambda)$ is
a $d$-th order twin-free subset of the set $\{1,2,\ldots, s-1\}$
by Lemma~\ref{lemma}. That is, $N_{d}(s)= \left| X_s \right| $.
Suppose that $X \in X_s$. If $s-1 \in X$, then $s-2, s-3, \ldots, s-(d+1) \notin X$,
since $X$ is a $d$-th order twin-free set. So,
$$ \left| \{ X \in X_s : (s-1) \in X \} \right| = \left| X_{s-(d+1)} \right|,$$ and
$$ \left| \{ X \in X_s : (s-1) \notin X \} \right| = \left| X_{s-1} \right| .$$
Thus, $\left| X_s \right|$ = $\left| X_{s-1} \right| +\left| X_{s-(d+1)} \right|$.
Notice that
\begin{eqnarray*}
N_{d}(1) & = & \left| X_1 \right| = 1 \\
N_{d}(2) & = & \left| X_2 \right|=2 \\
\vdots     & \vdots & \vdots  \\
N_{d}(d) & = & \left| X_d \right|=d \\
N_{d}(d+1)& = & \left| X_{d+1} \right|=d+1.
\end{eqnarray*}
So, we get the required result.
\end{prf}

If we take the value $d=1$ in Theorem~\ref{number}, we get
the number of $(s,s+1)$-core partitions with distinct parts is
the Fibonacci number $F_{s+1}$ in \cite{straub,xiong}.

\begin{example}
	For $d=2$, $N_{2}(6)=9$. The seven $(6,7)$-core partitions with
	$2$-distinct parts are
	$$\{ \}, \{1\},\{2\},\{3\}, \{3, 1\}, \{4\}, \{4, 1\}, \{5\}, \{4, 2\}.$$ 
	
	We can see from Table 1 that the number $N_{2}(s)$ of $(s,s+1)$-core partitions with $2$-distinct
	parts for $1 \le s \le 8$.
	
	\begin{table}[h]\caption{The number $N_{2}(s)$ of $(s,s+1)$-core
			partitions with $2$-distinct parts}
		\begin{center}
			\begin{tabular} {c||c|c|c|c|c|c|c|c} \label{table}
				$s$          & 1 & 2 & 3 & 4 & 5 & 6 & 7& 8      \\ \hline
				$N_{2}(s)$   & 1 & 2 & 3 & 4 & 6 & 9 & 13& 19 
			\end{tabular}
		\end{center}
	\end{table} 
	
	The generating function of the sequence $N_{2}(s)$ is
	
	$$-\dfrac{x^2+x+1}{x^3+x-1}.$$
	
	Also, the sequence $N_{2}(s)$ satisfies the recurrence relation
	
	$$N_{2}(s)=N_{2}(s-1)+N_{2}(s-3).$$
	
\end{example}

\begin{theorem}\label{size}
	If $s \equiv 0,1 \: \text{or} \: \: 2 \: (mod~{d+2})$ then
	the largest size of $(s,s+1)$-core partitions with
	$d$-distinct parts is
	$$\left[\dfrac{1}{d+2}\binom{s+1}{2}+\dfrac{s(d-1)}{2(d+2)}\right],$$
	otherwise
	$$\left[\dfrac{1}{d+2}\binom{s+1}{2}+\dfrac{s(d-1)}{2(d+2)}+1\right]$$
	where $[x]$ is the largest integer not greater than $x$.
	
\end{theorem}

\begin{prf}
	Let $\lambda$ is a $(s,s+1)$-core partitions with
	$d$-distinct parts. Suppose that $\beta(\lambda)=\{x_1,x_2,\ldots,x_k\}$.
	We need to maximize $\lambda$ and since $\beta(\lambda)$ is a $d$-th order twin-free set, we need $x_1=s-1, x_2=s-1-(d-1)$ and generally $x_i=s-d(i-1)-i$, so
	\begin{eqnarray*}
	|\lambda|&=& \sum_{i=1}^{k}x_i-\binom{k}{2} \\
	         &\le & \sum_{i=1}^{k}\left(s-d(i-1)-i\right)-\binom{k}{2} \\
	         &=& sk+\dfrac{dk-dk^2}{2}-k^2.
	\end{eqnarray*}
	
	Also, to maximize $\lambda$, we want to take $k$ as large as possible, however we also have to subtract the $\binom{k}{2}$ term. So if $x_k < (k-1)=\binom{k}{2}-\binom{k-1}{2}$, the gain we have made by include $x_k$ is offset by the loss of the second term.
	So, there are sometimes two $(s,s+1)$ cores with $d$-distinct parts and maximal size: this is when we have $x_k=k-1$, and so it makes no difference whether we include this term or not.
	
	When $s=(d+2)n$ for some integer $n$, we obtain
	$$|\lambda| \le sk+\dfrac{dk-dk^2}{2}-k^2 \le \dfrac{(d+2)n^2}{2}+\dfrac{dn}{2}.$$  
	When $s=(d+2)n+r,\: 1 \le r \le d+1$ for some integer $n$, we obtain
	$$|\lambda| \le sk+\dfrac{dk-dk^2}{2}-k^2 \le \dfrac{(d+2)n^2}{2}+\dfrac{dn}{2}+rn+(r-1).$$
	So, we can get the desired result for each case.  
\end{prf}

If we take the value $d=1$ in Theorem~\ref{size}, we get
the largest size of the $(s,s+1)$-core partitions with distinct
parts is $\left[\dfrac{1}{3}\binom{s+1}{2}\right]$ in \cite{xiong}.

\begin{example}
	For $s=6$ and $d=2$, since $s \equiv 2~(mod~4)$, the largest size of $(6,7)$-core
	partitions with $2$-distinct parts is
	$$\left[\dfrac{1}{2+2}\binom{6+1}{2}+\dfrac{6(2-1)}{2(2+2)}\right]=6$$
	by Theorem~\ref{size}. Indeed, $(6,7)$-core partitions with $2$-distinct parts are
	$$\{ \}, \{1\},\{2\},\{3\}, \{3, 1\}, \{4\}, \{4, 1\}, \{5\}, \{4, 2\}.$$
	So, the largest size of $(6,7)$-core partitions with $2$-distinct parts is $4+2=6$.
	
	For $s=7$ and $d=2$, since $s \equiv 3~(mod~4)$, the largest size of $(7,8)$-core
	partitions with $2$-distinct parts is
	$$\left[\dfrac{1}{2+2}\binom{7+1}{2}+\dfrac{7(2-1)}{2(2+2)}+1\right]=8$$
	by Theorem~\ref{size}. Indeed, $(7,8)$-core partitions with $2$-distinct parts are
	$$\{ \}, \{1\},\{2\},\{3\}, \{3, 1\}, \{4\}, \{4, 1\}, \{5\}, \{4, 2\}, \{5, 1\}, \{6\}, \{5, 2\}, \{5, 3\}.$$
	So, the largest size of $(7,8)$-core partitions with $2$-distinct parts is $5+3=8$.

\end{example}

\begin{theorem}\label{numbersize}
If $s \equiv 1\: (mod~{d+2})$ then there are two
$(s,s+1)$-core partitions of largest size with
$d$-distinct parts, otherwise there is only one
such partition of largest size.
\end{theorem}

\begin{prf}
Note that if $\lambda$ is an $(s,s+1)$-core partition with
$d$-distinct parts who has the largest size, then
$\beta(\lambda)=\{ s-1,s-(d+2),\ldots, s-((k-1)d+k)\}$
for some integer $k$. When $t=(d+2)n$ for some integer $n$,
$\lambda$ has the largest size if and only if $k=n$. When
$t=(d+2)n+1$ for some integer n, $\lambda$ has the largest size
if and only if $k=n$ or $k=n+1$. For all other cases $t=(d+2)n+r$,
where $2 \le r \le d+1$, $\lambda$ has the largest size
if and only if $k=n+1$. So, we get the desired result. 
\end{prf}

If we take the value $d=1$ in Theorem~\ref{numbersize}, we get
the number of the largest size of the $(s,s+1)$-core partitions with distinct
parts is $\dfrac{3-(-1)^{s \mod 3}}{2}$ in \cite{xiong}.

\begin{example}
For $s=5$ and $d=2$, since $s \equiv 1~(mod~4)$, there are only two
$(s,s+1)$-core partitions of largest size with $2$-distinct parts
by Theorem~\ref{numbersize}. Actually,
$(5,6)$-core partitions with $2$-distinct parts are
$$\{ \}, \{1\},\{2\},\{3\}, \{3, 1\}, \{4\}.$$
So, there are two partitions of the largest size of $(6,7)$-core partitions with $2$-distinct parts.
These partitions are $\{3, 1\}$ and $\{4\}$.

For $s=8$ and $d=3$, since $s \equiv 3~(mod~5)$, there are only one
$(s,s+1)$-core partition of the largest size with $3$-distinct parts
by Theorem~\ref{numbersize}. Indeed, $(8,9)$-core partitions with $3$-distinct parts are
$$\{ \}, \{1\},\{2\},\{3\}, \{4\}, \{4, 1\}, \{5\}, \{5, 1\}, \{6\}, \{5, 2\}, \{6, 1\}, \{7\}, \{6, 2\}, \{6, 3\}.$$
So, there are only one partition of the largest size of $(8,9)$-core partitions with $3$-distinct parts.
This partition is $\{6, 3\}$.
\end{example}

\section{$(s,s+r)$-core partitions with $d$-distinct parts}

More generally, we propose a conjecture about the number of
$(s,s+r)$-core partitions with $d$-distinct parts for $1 \le r \le d$.
This conjecture is based on experimental evidence and has been verified
for $s < 10$ after listing all relevant partitions. We will present
some of our experimental results in Table~\ref{table1} and Table~\ref{table2}.

Table~\ref{table1} shows $(s,s+2)$-core partitions
with $d$-distinct partitions for $ 2 \leq d \leq 7$.
  
\vspace*{0.5 cm}  
\begin{table}[h]
	\caption{The number of $(s,s+2)$-core partitions with $d$-distinct parts}\label{table1}
	\centering
	\begin{tabular}{ *{9}{|c}|} 
		\hline
		\backslashbox{$d$}{$(s,s+2)$} & $(1,3)$ & (2,4) & (3,5) & (4,6) & (5,7) & (6,8) & (7,9) & (8,10)\\
		\hline
		2 & 1  & 2  & 4 & 5 & 7 & 11 & 16 & 23 \\
		3 & 1  & 2  & 3 & 5 & 6 & 8 & 11 & 16 \\
        4 &  1  & 2  & 3 & 4 & 6 & 7 & 9 & 12 \\
        5 &  1  & 2  & 3 & 4 & 5 & 7 & 8 & 10 \\
        6 &  1  & 2  & 3 & 4 & 5 & 6 & 8 & 9 \\
        7 &  1  & 2  & 3 & 4 & 5 & 6 & 7 & 9\\
		\hline
	\end{tabular}
\end{table}

\vspace*{0,5 cm}

Table~\ref{table2} shows $(s,s+3)$-core partitions
with $d$-distinct partitions for $ 3 \leq d \leq 7$.  
According to our experiments, we present the following conjecture.
\vspace*{0.5 cm}

\begin{table}[h]
	\caption{$(s,s+3)$-core partitions with $d$-distinct parts}\label{table2}
	\centering
	\begin{tabular}{ *{9}{|c}|} 
		\hline
		\backslashbox{$d$}{$(s,s+3)$} & (1,4) & (2,5) & (3,6) & (4,7) & (5,8) & (6,9) & (7,10) & (8,11)\\
		\hline
		3 & 1  & 2  & 3 & 6 & 7 & 9 & 12 & 18 \\
		4 &  1  & 2  & 3 & 4 & 7 & 8 & 10 & 13 \\
		5 &  1  & 2  & 3 & 4 & 5 & 8 & 9 & 11 \\
		6 &  1  & 2  & 3 & 4 & 5 & 6 & 9 & 10 \\
		7 &  1  & 2  & 3 & 4 & 5 & 6 & 7 & 10 \\
		\hline
	\end{tabular}
\end{table}

\begin{conjecture}
	For $1 \le r \le d$, the number $N_{d,r}(s)$ of $(s,s+r)$-core partitions with
	$d$-distinct parts is characterized by
	$N_{d,r}(s)=s$ for $1 \le s \le d$, $N_{d,r}(d+1)=d+r$, and for $s \ge d+2$,
	$$ N_{d,r}(s)=N_{d,r}(s-1)+N_{d,r}(s-(d+1)).$$
\end{conjecture}

\begin{example}
For $s=6, d=3$ and $r=2$, the eight $(s,s+r)$-core, i.e. $(6,8)$-core, partitions with
$3$-distinct parts are
$$\{ \}, \{1\},\{2\},\{3\},\{4\},\{5\},\{1,4\},\{1,6\}.$$ 

We can see from Table 1 that the number $N_{3,2}(s)$ of $(s,s+2)$-core partitions with $3$-distinct
parts for $1 \le s \le 9$.

\begin{table}[h]\caption{The number $N_{3,2}(s)$ of $(s,s+2)$-core
partitions with $3$-distinct parts}
\begin{center}
\begin{tabular} {c||c|c|c|c|c|c|c|c|c} \label{table}
$s$          & 1 & 2 & 3 & 4 & 5 & 6 & 7& 8 & 9     \\ \hline
$N_{3,2}(s)$ & 1 & 2 & 3 & 5 & 6 & 8 & 11& 16 & 22
\end{tabular}
\end{center}
\end{table} 

The generating function of the sequence $N_{3,2}(s)$ is

$$-\dfrac{2x^3+x^2+x+1}{x^4+x-1}.$$

Also, the sequence $N_{3,2}(s)$ satisfies the recurrence relation

$$N_{3,2}(s)=N_{3,2}(s-1)+N_{3,2}(s-4).$$

\end{example}

\hfill \rule{0.5em}{0.5em}

\end{document}